\documentclass[a4paper,leqno,11pt]{article}
\usepackage{amssymb,amsfonts,amsmath,amsthm}
\usepackage{epsfig,epsf,curves}
\usepackage{bbm,epic,eepic}
\usepackage[mathscr]{eucal}
\theoremstyle{remark}

\begin{document}

\sloppy


\def\Mat#1#2#3#4{\left(\!\!\!\begin{array}{cc}{#1}&{#2}\\{#3}&{#4}\\ \end{array}\!\!\!\right)}

\def\tAA{{\Bbb A}}
\def\CC{{\Bbb C}}
\def\HH{{\Bbb H}}
\def\NN{{\Bbb N}}
\def\QQ{{\Bbb Q}}
\def\RR{{\Bbb R}}
\def\ZZ{{\Bbb Z}}
\def\FF{{\Bbb F}}
\def\SS{{\Bbb S}}
\def\GG{{\Bbb G}}
\def\PP{{\Bbb P}}
\def\LL{{\Bbb L}}


\title{Conjugation of Hilbert modular forms and trace formula}

\author{Joachim Mahnkopf}

\maketitle

{\bf Abstract} {We describe (in a representation theoretic setting) a simple comparison of trace formulas, which implies that the conjugate of a Hilbert modular form $f$ by an automorphism of ${\Bbb C}$ again is a Hilbert modular form of the same level and conjugate weight as $f$. This is a Theorem of Shimura for which we obtain a new proof (cf. Theorem 3.3 and Corollary 3.4).

}

\vspace{1.5cm}

\centerline{\bf \large Introduction}

{\bf 0.1.} Let $f$ denote a Hilbert modular form of even weight ${\bf k}$. The conjugate ${^\sigma}f$ of $f$ by an automorphism $\sigma$ of ${\Bbb C}$ is obtained by applying $\sigma$ to the coefficients of the Fourier expansion of $f$. Thus, ${^\sigma}f$ is a holomorphic mapping and Shimura's Theorem states that  ${^\sigma}f$ again is a Hilbert modular form of conjugate weight ${^\sigma}{\bf k}$ (cf. [Sh]). A different proof of
Shimura's Theorem has been given by Garrett and is described in his book [G]. We note that in the generality of ${\rm GL}_n$ Clozel has proven a representation theoretic analogoue of Shimura's Theorem for algebraic automorphic representations, which satisfy a certain regularity condition (cf. [C]).

In this article we describe in the case of ${\rm GL}_2$ over totally real fields a comparison of trace formulas, which implies that conjugation by automorphisms of ${\Bbb C}$ preserves the automorphic property. We use a representation theoretic formulation. To be more precise, we let $F/{\Bbb Q}$ be a totally real extension with adele ring ${\Bbb A}$. We denote by $L_{\rm cusp}^2({\bf k})$ the space of (adelic) Hilbert cusp forms on ${\rm GL}_2({\Bbb A})$ of weight ${\bf k}=(k_v)_{v\not|\infty}$. $L_{\rm cusp}^2({\bf k})$ is a module under the Hecke algebra ${\cal H}_f$ of ${\rm GL}_2({\Bbb A}_f)$ and we denote by $\Pi({\bf k})$ is the set of all irreducible representations of ${\cal H}_f$, which appear in $L_{\rm cusp}^2({\bf k})$ (in the main part we will also fix a nebentype $\omega$ and a level $K$). Now, if the $\sigma$-conjugate of any Hilbert cusp form again is a Hilbert cusp form of conjugate weight ${^\sigma}{\bf k}=(k_{\sigma^{-1}(v)})_{v\not|\infty}$ then conjugation of abstract representations $\pi\mapsto{^\sigma}\pi$ (cf. section 2.1 for the definition) defines a map 
$$
\sigma:\,\Pi({\bf k})\rightarrow\Pi({^\sigma}{\bf k}).\leqno(0.1)
$$ 
Explicitly, this means that if $\pi$ is the finite part of a cuspidal representation, whose archimedean component has $K_\infty$-type given by ${\bf k}$ then ${^\sigma}\pi$ is the finite part of a cuspidal representation whose archimedean component has $K_\infty$-type given by ${^\sigma}{\bf k}$. Our aim is to establish the existence of the map (0.1) by a comparison of trace formulas. To this end we define a corresponding (dual) map on the 
Hecke algebra 
$$
\sigma:\,{\cal H}_f\rightarrow{\cal H}_f
$$
by simply setting ${^\sigma\varphi}(x)=\sigma(\varphi(x))$. The following result then connects conjugation of automorphic representations to the trace formula.

\medskip

{\bf Proposition 0.1. } (cf. Proposition 2.1) {\it Let $\sigma\in{\rm Aut}({\Bbb C}/{\Bbb Q})$. If for any $\varphi\in{\cal H}_f$ 
$$
\sigma\,{\rm tr}\,\varphi|_{L_{\rm cusp}^2({\bf k})}={\rm tr}\,{^\sigma \varphi}|_{L_{\rm cusp}^2({^\sigma}{\bf k})}\leqno(0.2)
$$
then conjugation defines a map 
$$
\sigma:\,\Pi({\bf k})\rightarrow\Pi({^\sigma {\bf k}}).
$$
}

\medskip

Comparing trace formulas for $\varphi$ in weight ${\bf k}$ and for ${^\sigma}\varphi$ in weight $^\sigma{\bf k}$ we verify under a certain {\it algebraicity condition} that the trace identity (0.2) holds. Thus, we obtain a proof of Shimura's Theorem (cf. Theorem 3.3 and Corollary 3.4). 

\medskip

{\bf 0.2. } We add some Remarks. 1.) The trace formula is some kind of universal principle for proving the existence of maps between sets of automorphic representations. In this sense it is natural to try to prove Shimura's Theorem by establishing the existence of (0.1) via a comparison of trace formulas.
2.) The proofs of Shimura and Clozel make use of a ${\Bbb Q}$-structure on a certain space which contains the Hecke module of automorphic cusp forms: Shimura uses the 
${\Bbb Q}$-structure given by the $q$-expansion, Clozel uses the ${\Bbb Q}$-structure on deRham cohomology given by singular cohomology via the deRham isomorphism. In our proof we use a different ${\Bbb Q}$-structure: it is defined on the space of all mappings from ${\rm GL}_2({\Bbb A}_f)$ to ${\Bbb C}$ by the subspace of mappings, which are ${\Bbb Q}$-valued (cf. equation (0.2)). 
3.) The proof of equation (0.2) reduces to showing the algebraicity of the local archimedean distributions appearing on the geometric side of trace formula (cf. Lemma 2.2 and Corollary 3.6). In the case of ${\rm GL}_2$ this would have been possible using explicit calculations. Instead we will use general principles from harmonic analysis based on [Ca]. We hope that in this way the proof will generalize to higher rank groups as ${\rm Sp}_n$.  The trace formula may be seen here as a device which converts the local considerations about archimedean orbital integrals into global existence statements about Hilbert cusp forms. We note that the Multiplicity $1$ Theorem enables to use a simple trace formula. Nevertheless it may be interesting to examine the algebraicity of all the distributions appearing in geometric side of the Selberg trace formula. We discuss this briefly in section 3.3

\medskip

{\bf 0.3.} In section 1, for convenience, we review some well known facts about Hilbert modular forms and representations of ${\rm GL}_2$. In section 2 we prove Proposition 0.1. In fact we prove a slightly stronger result
(cf. Proposition 2.1), which will enable us in section 3 to verify (0.2) by using a simple trace formula.

\section {Hilbert modular forms. }

{\bf 1.1. Notations. } We fix a totally real number field $F/{\Bbb Q}$ with Galois group ${\cal G}={\rm Gal}(F/{\Bbb Q})$. We denote by ${\cal O}$
the integers of $F$ and by ${\Bbb A}$ its ring of adeles. By a place $v$ of $F$ we understand an equivalence class of valuations
$|\cdot|:\,F\rightarrow {\Bbb R}$ and we denote by $S_\infty$ the set of archimedean places of $F$ and by $i_v:\,F\hookrightarrow F_v$ the completion of $F$ at the place 
$v$. ${\cal G}$ acts on the places of $F$: if $v$ is represented by the valuation $|\cdot|$ then
$\tau v$ is the place, which is represented by $|\cdot|\circ\tau$, $\tau\in{\cal G}$. We fix an archimedean place $v_0$; the elements in $S_\infty$ then 
are given as $\tau v_0$, $\tau\in{\cal G}$, where the $\tau v_0$ are pairwise different. 
We note that $\tau\in{\cal G}$ extends to a morphism $\tau:\,F_v\rightarrow F_{\tau^{-1} v}$. 
%
%

We denote by $j$ the isomorphism $j:\,F_{v_0}\cong{\Bbb R}$; 
for any $\tau\in{\cal G}$ we obtain a commutative diagram
$$
\begin{array}{ccc}
{\Bbb R}&&\\
j\,\uparrow&&\\
F_{v_0}&\stackrel{\tau^{-1}}{\rightarrow}&F_{\tau v_0}\\
i_{v_0}\,\uparrow&&\uparrow\,i_{\tau v_0}\\
F&\stackrel{\tau^{-1}}{\rightarrow}&F.\\
\end{array}
$$

We set $G={\rm GL}_2/F$ and we denote by $Z/F$ the center of $G$. For any place $v$ we set $G_v=G(F_v)$ and $Z_v=Z(F_v)$ and if $v$ is
archimedean we denote by $Z_v^0$ the connected component of $1$ of $Z_v$. We set $G_\infty=\prod_{v\in S_\infty} G_v$ and 
$Z_ \infty=\prod_{v\in S_\infty} Z_v$. Similarly, we denote by $Z_2$ the center of ${\rm GL}_2({\Bbb R})$, by $Z_2^0$ its connected component 
containing $1$.

The morphisms $j,i_v$ and $\tau\in{\cal G}$ extend to $G(F_{v_0})$, $G(F)$ and $G(F_v)$ by applying them to the entries of a matrix and 
we obtain a diagram
$$
\begin{array}{ccccc}
&{\rm GL}_2({\Bbb R})&&&\\
j&\uparrow&&&\\
&G(F_{v_0})&\stackrel{\tau^{-1}}{\rightarrow}&G(F_{\tau v_0})&\\
i_{v_0}&\uparrow&&\uparrow&i_{\tau v_0}\\
&G(F)&\stackrel{\tau^{-1}}{\rightarrow}&G(F).&\\
\end{array}\leqno(1.0)
$$
We will use the following {\it identifications}: we identify any $\gamma\in G(F)$ with its image $ji_{v_0}(\gamma)$ in ${\rm GL}_2({\Bbb R})$ 
and we identify $G(F_{\tau v_0})$ with ${\rm GL}_2({\Bbb R})$ via the map $j\circ\tau$. The commutativity of (1.0) shows that under these 
identifications 
$$
i_{\tau v_0}(\gamma)=\tau(\gamma)\in {\rm GL}_2({\Bbb R})\leqno(1.1)
$$ 
for all $\gamma\in G(F)$. We identify $G(F_v)/Z_v^0={\rm GL}_2({\Bbb R})/Z^0_2={\rm SL}_2({\Bbb R})^\pm$ (matrices of determinant $\pm 1$) by 
sending $\gamma\in {\rm SL}_2({\Bbb R})^\pm$ to $\gamma Z_2^0\in{\rm GL}_2({\Bbb R})/Z_2^0$. 
Finally, under these identifications, to any function $\varphi:\,{\rm GL}_2({\Bbb R})\rightarrow{\Bbb C}$ and any place $v=\tau v_0$ 
corresponds the function 
$$
\hat{\varphi}_v=\varphi\circ j\circ\tau:\,G(F_{v})\rightarrow{\Bbb C}.\leqno(1.2)
$$

We fix a character $\omega:\,F^*\backslash {\Bbb A}^*\rightarrow {\Bbb C}^*$ and a compact open subgroup $K=\prod_{v\not\in S_\infty} K_v\le
\prod_{v\not\in S_|infty }G({\cal O}_v)$. For any finite place $v$ we denote by ${\cal H}_v(\omega_v,K_v)$ the local Hecke algebra consisting of 
all $K_v$-bi-invariant functions, which are compactly supported modulo center and have central character $\omega_v^{-1}$, 
i.e. $\varphi(zx)=\omega_v^{-1}(z)\varphi(x)$ for all $z\in Z(F_v)$. Thus, we have to assume that $Z(F_v)\cap K_v\le {\rm ker}\,\omega_v$,
which we can always achieve by intersecting $K_v$ with a sufficiently small principal congruence subgroup. At archimedean places we denote by ${\cal H}_v(\omega_v)$ the set of 
all compactly supported modulo center, smooth functions on $G(F_v)$ having central character $\omega_v^{-1}$. We also set 
${\cal H}_f(\omega_f,K)=\otimes_{v\not|\infty} {\cal H}_v(\omega_v,K_v)$ and  
${\cal H}_\infty(\omega_\infty)=\otimes_{v|\infty} {\cal H}_v(\omega_v)$. We further denote by 
${\cal H}({\rm GL}_2({\Bbb R}),\omega_{\Bbb R})$ the Hecke algebra consisting of functions, which are smooth and compactly supported modulo 
center and have central character $\omega_{\Bbb R}^{-1}$; we use a similar notation for the group ${\rm SL}_2({\Bbb R})^\pm$. If $(\pi,V)$ is 
any representation 
of $G({\Bbb A}_f)$ with central character $\omega_f$ we obtain a corresponding representation $(\pi,V^K)$ of the Hecke algebra 
${\cal H}_f(\omega_f,K)$ by setting
$$
\pi(\varphi)=\int_{G({\Bbb A}_f)/Z({\Bbb A}_f)} \varphi(g)\,\pi(g)\,dg
$$ 
for any $\varphi\in {\cal H}_f(\omega_f,K)$. An analogous statement holds for representations of 
$G(F_v)$ having central character $\omega_v$. We use the following notation: if $\gamma$ is contained
in a group $G$, then  we denote by $G(\gamma)$ the centralizer of $\gamma$ in $G$.

\medskip

{\it Remark. }{ Let $v=\tau v_0$, $\tau\in{\cal G}$, be any archimedean place. Using the above identifications we have for any compactly supported function $\varphi$
on ${\rm GL}_2({\Bbb R})$ and any $\gamma\in G(F)$
$$
\int_{{\rm GL}_2({\Bbb R})(\tau(\gamma))\backslash {\rm GL}_2({\Bbb R})} \varphi(x^{-1}\tau(\gamma)x)\,dx
=\int_{G(F_v)(i_{v}(\gamma))\backslash G(F_v)} \hat{\varphi}_v(x^{-1}i_{v}(\gamma)x)\,dx
$$
}

{\it Proof. } We calculate using the map $j\circ \tau$ and the commutativity of (1.0) in the last step
\begin{eqnarray*}
\int_{{\rm GL}_2({\Bbb R})(ji_{v_0}(\gamma))\backslash{\rm GL}_2({\Bbb R})} \varphi(x^{-1}ji_{v_0}(\gamma)x)\,dx
&=&\int_{G(F_v)(\tau^{-1} i_{v_0}(\gamma))\backslash G(F_v)} \varphi(j\tau(x)^{-1}ji_{v_0}(\gamma)j\tau(x))\,dx\\
&=&\int_{G(F_v)(\tau^{-1} i_{v_0}(\gamma))\backslash G(F_v)} \varphi\circ j\circ\tau(x^{-1}\tau^{-1} i_{v_0}(\gamma)x)\,dx\\
&=&\int_{G(F_v)(i_{\tau v_0}\tau^{-1}(\gamma))\backslash G(F_v)} \hat{\varphi}_v(x^{-1}i_{\tau v_0}\tau^{-1}(\gamma)x)\,dx.\\
\end{eqnarray*}
Taking into account our identification of $\gamma\in G(F)$ with its image $j i_{v_0}(\gamma)\le {\rm GL}_2({\Bbb R})$
and replacing $\tau^{-1}(\gamma)\in G(F)$ by $\gamma\in G(F)$ we obtain the claim of the Remark.

\medskip

{\bf 1.2. Projectors in the Hecke algebra. } We denote by $\delta_n:\,{\rm SO}_2({\Bbb R})\rightarrow{\Bbb C}^*$ the character which sends
$$
{\bf e}(\Theta)=\Mat{\cos\Theta}{\sin \Theta}{-\sin \Theta}{\cos \Theta}\mapsto e^{i n\Theta}.
$$
We let $(D_k,W_k)$ be the irreducible discrete series representation of ${\rm SL}_2({\Bbb R})$ of lowest weight $\delta_k$, $k\ge 2$; thus, 
$D_k$ has central character ${\rm sgn}^k:\,\Mat{-1}{}{}{-1}\rightarrow (-1)^k$. We let $(L_k,V_k)$ be the irreducible algebraic representation 
of ${\rm SL}_2$ of highest weight 
$$
\Mat{x}{}{}{x^{-1}}\rightarrow x^k,
$$ 
$k\ge 0$, with respect to the torus $T_2$ of diagonal matrices in ${\rm SL}_2$, i.e. $L_k$ has highest weight $k\gamma$, where $\gamma$ is the 
fundamental weight of ${\mathfrak sl}_2$ corresponding to the Borel subalgebra of upper triangular matrices. $L_k$
induces representations of ${\rm SL}_2({\Bbb C})$ and ${\rm SL}_2({\Bbb R})$. Moreover, since $L_k^{\Mat{-1}{}{}{1}}\cong L_{k}$
the representation $L_k$ extends to a representation of ${\rm SL}_2({\Bbb R})^\pm$.

\medskip

{\bf Lemma 1.1. }{\it There is $\varphi_k\in {\cal H}({\rm sgn}^k,{\rm SL}_2({\Bbb R})^\pm)$ such that
$$
{\rm tr}\,\pi(\varphi_k)={\rm tr}\,\int_{{\rm SL}_2({\Bbb R})^\pm/\{\pm 1\}}\varphi_k(g)\,\pi(g)\, dg=0
$$
if $\pi$ is an irreducible representation of the principal series, or $\pi\cong D_n$ with $n\not=k$ or $\pi\cong L_n$ with $n\not=k-2$ and
$$
{\rm tr}\,\pi(\varphi_k)=\left\{\begin{array}{ccc}
1&{\rm if}&\pi\cong D_k\\
-1&{\rm if}&\pi\cong L_{k-2}\\
\end{array}\right.
$$
More precisely, $D_k(\varphi_k)$ leaves the $\delta_n$-isotypical components $W_k(\delta_n)\le W_k$ invariant and
$$
{\rm tr}\,D_k(\varphi_k)|_{W_k(\delta_n)}=\left\{\begin{array}{ccc}
1&{\rm if}&n=k\\
0&{\rm if}&k>n.\\
\end{array}\right.
$$

}

{\it Proof. } We follow the argument given in [Ca], p. 149/150 in the case $k=2$ and $G={\rm PSL}_2$, which we extend to the case $k\ge 2$ and ${\rm GL}_2$. We first let $k\ge 0$. Cartan
decomposition implies that there is a compactly supported 
function $f_k:\,{\rm SL}_2({\Bbb R})\rightarrow {\Bbb C}$ such that 
$$
f_k({\bf e}(\Theta_1)g{\bf e}(\Theta_2))=\delta_{k}(\Theta_1)\delta_{k}(\Theta_2)\,f_k(g)
$$ 
for all $g\in{\rm SL}_2({\Bbb R})$. We note that this implies that
$$
f_k(\Mat{-1}{}{}{-1})=(-1)^k,
$$
i.e. $f_k\in {\cal H}({\rm sgn}^k,{\rm SL}_2({\Bbb R}))$. $f_k$ vanishes on all $\delta_n$-isotypical components with $n\not=k$ and
leaves the $\delta_k$-isotypical component of any representation invariant. Moreover, if we choose the support of $f_k$ sufficiently close to 
${\rm SO}_2({\Bbb R})$ and suitably normalize $f_k$, we obtain
$$
D_k(f_k)|_{W_k(\delta_k)}={\rm id}.
$$
We set $g_k=f_k-f_{k-2}$, $k\ge 2$. Since $D_n=\bigoplus_{m} {D_n(\delta_m)}$, where $m=\pm n,\pm (n+2),\pm (n+4),\ldots$ we obtain
$$
{\rm tr}\,\pi(g_k)=\left\{\begin{array}{ccc}
0&{\rm if}&\pi\cong D_n\;\mbox{with}\;n\not=k\\
1&{\rm if}&\pi\cong D_k.\\
\end{array}\right.\leqno(1.3)
$$
More precisely, we obtain
$$
{\rm tr}\,D_k(g_k)|_{W_k(\delta_n)}=\left\{\begin{array}{ccc}
0&{\rm if}&n>k\\
1&{\rm if}&n=k.\\
\end{array}\right.\leqno(1.3')
$$
As in [loc. cit.] we choose a compactly supported ${\rm SO}_2({\Bbb R})$-bi-invariant function $h_k$ on ${\rm SL}_2({\Bbb R})$ such that the 
Harish-Chandra transform $H_{h_k}$ of $h_k$ equals the Harish-Chandra transform $H_{g_k}$ of $g_k$. We set $\varphi_k=g_k-h_k$, hence 
$H_{\varphi_k}=0$. Since no discrete series representation has trivial ${\rm SO}_2({\Bbb R})$-types, equation (1.3) and (1.3') remain valid 
for $\varphi_k$ (note that $h_k$ is ${\rm SO}_2({\Bbb R})$ bi-invariant). On the other hand for any principal series representation $\pi_\chi$ 
of ${\rm SL}_2({\Bbb R})$ we have
$$
{\rm tr}\,\pi_\chi(\varphi_k)=\int_{T_2({\Bbb R})} H_{\varphi_k}(t)\chi(t)\,dt=0.\leqno(1.4)
$$
The exact sequence
$$
0\rightarrow L_{n-2}\rightarrow{\rm Ind}(|\cdot|_\infty^{(n-1)/2},|\cdot|_\infty^{-(n-1)/2}{\rm sgn}^n)\rightarrow D_{n}\rightarrow 0
$$
together with the vanishing of ${\rm tr}\,\pi(\varphi_k)$ for representations $\pi$ of the principal series and equation (1.3) shows that
$$
{\rm tr}\,L_n(\varphi_k)=\left\{\begin{array}{ccc}
-1&{\rm if}&n=k-2\\
0&&{\rm else.}\\ 
\end{array}\right.\leqno(1.5)
$$
Thus, altogether, $\varphi_k$ satisfies the requirements of the Lemma, except that it is contained in the Hecke algebra of ${\rm SL}_2({\Bbb R})$. 
We extend $\varphi_k$ to a function on ${\rm SL}_2({\Bbb R})^\pm$ by setting it equal to $0$ on the connected component 
${\rm SL}_2({\Bbb R})^-$. Equations (1.3), (1.3'), (1.4), (1.5), then imply that $\varphi_k$ satisfies the requested properties. Thus, the proof of the Lemma is complete.

\medskip

We extend Lemma 1.1 to ${\rm GL}_2$. This will involve an algebraicity condition. For any pair of integers $k,w\in{\Bbb Z}$ with $k\ge 2$ we denote by $(D_{k,w},W_{k,w})$ the 
irreducible discrete series representation of ${\rm GL}_2({\Bbb R})$, which has lowest ${\rm SO}_2({\Bbb R})$-type $\delta_{k}$ and whose 
central character $\omega_{D_{k,w}}$ when restricted to the connected component $Z_2^0$ is given by $x^w:\,\Mat{x}{}{}{x}\mapsto x^w$. Thus, $D_{k,w}$
is a quotient of ${\rm Ind}(|\cdot|_\infty^{(k-1+w)/2},|\cdot|_\infty^{(-k+1+w)/2}{\rm sgn}^k)$. We note that the 
central character satisfies 
$$
\omega_{D_{k,w}}(-1)=(-1)^k,\leqno(1.6)
$$ 
hence, $\omega_{D_{k,w}}$ and therefore $D_{k,w}$ is uniquely determined by $k$ and $w$ and any discrete series representation is isomorphic
to a representation $D_{k,w}$. We assume that the following algebraicity condition holds
$$
k\equiv w\pmod{2}\leqno(Alg_{\infty})
$$ 
This has the following two consequences:

1. $D_{k,w}$ has {\it algebraic} central character 
$$
x^w:\,\Mat{x}{}{}{x}\mapsto x^w\qquad (x\in{\Bbb R}^*).\leqno(1.7) 
$$

2. $D_{k,w}$ fits in an exact sequence
$$
0\rightarrow L_{k-2,w}\rightarrow{\rm Ind}(|\cdot|_\infty^{(k-1+w)/2},|\cdot|_\infty^{(-k+1+w)/2}{\rm sgn}^k)\rightarrow D_{k,w}\rightarrow 0.
$$  
where now $L_{k-2,w}$ is an {\it algebraic} representation of ${\rm GL}_2$. More precisely, $(L_{k,w},V_{k,w})$ is the finite dimensional, 
algebraic irreducible representation of ${\rm GL}_2$ of highest weight $\Mat{x}{}{}{y}\mapsto x^{(k+w)/2}y^{(-k+w)/2}$.


\medskip

We denote by $\varphi_{k,w}$ the image of $\varphi_k$ under the canonical isomorphism ${\cal H}({\rm SL}_2({\Bbb R})^\pm,{\rm sgn}^k)
\cong {\cal H}({\rm GL}_2({\Bbb R}),x^w)$.

\medskip

{\bf Corollary 1.2. } {\it Assume that ($Alg_\infty)$ holds. Let $\pi$ be an irreducible representation of ${\rm GL}_2({\Bbb R})$, whose central 
character satisfies $\omega_\pi(\Mat{x}{}{}{x})=x^w$. Then
$$
{\rm tr}\,\pi(\varphi_{k,w})={\rm tr}\,\int_{{\rm GL}_2({\Bbb R})/Z_2}\varphi_{k,w}(g)\,\pi(g)\, dg=0
$$
if $\pi$ is an irreducible representation of the principal series, or $\pi\cong D_{n,w}$ with $n\not=k$ or $\pi\cong L_{n,w}$ with $n\not=k-2$ 
and
$$
{\rm tr}\,\pi(\varphi_{k,w})=\left\{\begin{array}{ccc}
1&{\rm if}&\pi\cong D_{k,w}\\
-1&{\rm if}&\pi\cong L_{k-2,w}.\\
\end{array}\right.
$$
More precisely, we have
$$
{\rm tr}\,D_{k,w}(\varphi_{k,w})|_{W_{k,w}(\delta_n)}=\left\{\begin{array}{ccc}
1&{\rm if}&k=n\\
0&{\rm if}&n>k.\\
\end{array}\right.
$$

}

The Corollary is an immediate consequence of Lemma 1.1. We note that the integration is well defined because ($Alg_\infty$) implies that
$\varphi_{k_v,w}$ has central character $x^{-w}$.

{\bf 1.3. Hilbert modular forms. } We fix a compact subgroup $K=\prod_{v\not\in S_\infty} K_v\le\prod_{v\not\in S_\infty}G({\cal O}_v)$, a 
character $\omega:\,F^*\backslash {\Bbb A}^*\rightarrow{\Bbb C}^*$ and weight ${\bf k}=(k_v)_{v\in S_\infty}$, $k_v\in{\Bbb N}$, $k_v\ge 2$. 
Since $\omega$ is an idele class character there is an integer $w\in{\Bbb Z}$ such that 
$$
\omega_v(x)=\pm |x|_v^w\leqno(1.8)
$$ 
for all $v\in S_\infty$ and $x\in F_v^*$. We denote by $L_d^2(\omega)$ the discrete part of the spectrum of the $G({\Bbb A})$-module $L^2(G(F)\backslash G({\Bbb A}),\omega)$, which 
consists of square integrable functions having central character $\omega$ and $L^2_0(\omega)\le L^2_d(\omega)$ is the subspace of cusp 
forms. We set $D_{{\bf k},w}=\otimes_{v\in S_\infty} D_{k_v,w}$, where we view $D_{k_v,w}$ as a representation of $G(F_v)$ via the
identifications in section 1.1, and $\delta_{\bf k}=\otimes_{v\in S_\infty} \delta_{k_v}$; thus, $D_{{\bf k},w}$ resp. $\delta_{\bf k}$ is a representation of $G_\infty$ resp. of  ${\rm SO}_{2,\infty}=\prod_{v\in S_\infty} {\rm SO}_2({\Bbb R})$. We then define 
$$
L_0^2(\omega,{\bf k})=\bigoplus_{(\pi,V_\pi)\le L_0^2(\omega)\atop \pi_\infty\cong D_{{\bf k},w}} V_\pi(\delta_{\bf k}).
$$
Thus, $L_0^2(\omega,{\bf k})$ is the space of adelic Hilbert modular forms of weight ${\bf k}$ and central character
$\omega$.  
$L_0^2(\omega,{\bf k})$ is a $G({\Bbb A}_f)$-module and the subspace of $K$-invariant vectors $L_0^2(\omega,{\bf k})^K$ is a 
${\cal H}_f(\omega_f,K)$-module.

\medskip

{\it Remark. } Since $D_{k_v,w}$ has central character satisfying $\omega_{D_{k_v,w}}(-1)=(-1)^{k_v}$ (cf. (1.6)) we see that 
$L_0^2(\omega,{\bf k})^K$ is the empty sum unless 
$$
\omega_v(-1)=(-1)^{k_v}\leqno(1.9)
$$ 
for all $v\in S_\infty$. From now on we shall always assume that this holds.


\medskip

{\bf 1.4. Algebraic Hilbert modular forms. } From now on we shall always assume that the weight ${\bf k}$ and the character $\omega$ 
satisfy the following algebraicity condition: 
$$
k_v\equiv w\pmod{2}\quad \mbox{for all}\; v\in S_\infty.\leqno(Alg)
$$
(cf. equation (Alg$_\infty$)). We note that condition (Alg) is the algebraicity condition of [C], Definition 1.8 in the special case of 
${\rm GL}_2/F$. Condition (Alg) together with equations (1.8) and (1.9) imply that 
$$
\omega_v(x)=x^w\leqno(1.10)
$$ 
for all $v\in S_\infty$, $x\in F_v^*$. In particular, $\omega$ is an algebraic character. 
Finally, we denote by 
$$
\Pi_{\bf k}(\omega,K)
$$ 
the set of all representations $\pi_f$ of $G({\Bbb A}_f)$ such that $\pi_f^K$ appears in $L_0^2(\omega,{\bf k})^K$. I.e. $\Pi_{\bf k}(\omega,K)$
consists of the finite parts of cuspidal automorphic representations of $G({\Bbb A})$, whose infinity component has ${\rm SO}_2({\Bbb R})$-type 
$\delta_{\bf k}$. Using equation (1.2) we define the following element (cf. equation (1.10)): 
$$
\hat{\varphi}_{{\bf k},w}=\otimes_{v\in S_\infty}\hat{\varphi}_{k_v,w,v}\in{\cal H}_\infty(\omega_\infty).
$$

\medskip

{\bf Corollary 1.3. } {\it For all $\varphi\in{\cal H}_f(K)$ we have 
$$
{\rm tr}\, \varphi|_{L_0^2(\omega,{\bf k})}={\rm tr}\,\varphi\otimes\hat{\varphi}_{{\bf k},w}|_{L_0^2(\omega)}.
$$

}

\section{Conjugation and Trace.}

{\bf 2.1. Conjugation of representations. } For the moment $(\pi,W)$ denotes any representation of $G({\Bbb A}_f)$ or of $G(F_v)$. For any $\sigma\in{\rm Aut}({\Bbb C}/{\Bbb Q})$ we define 
the conjugate representation $({^\sigma}\pi,{^\sigma W})$ as follows. First, we define the ${\Bbb C}$-vector space ${^\sigma}W$: we 
set ${^\sigma W}=W$ as abelian groups and using the scalar multiplication "$\cdot$" on $W$ we define a scalar multiplication "$\cdot_\sigma$" 
on ${^\sigma W}$ as $\alpha\cdot_\sigma v=\sigma^{-1}(\alpha)v$ for all $\alpha\in{\Bbb C}$ and $v\in {^\sigma}W$. We then define a representation 
$^\sigma\pi$ on $^\sigma W$ by setting ${^\sigma}\pi(g)(v)=\pi(g)(v)$ for all $g\in G({\Bbb A}_f)$ or $g\in G(F_v)$ and $v\in{^\sigma W}$. The corresponding 
representation of the Hecke algebra is given by ${^\sigma}\pi(\varphi)v=\int \sigma^{-1}(\varphi(g))\, \pi(g)v\,dg$. 

In this section we want to show that conjugation of representations preserves the automorphic property; more precisely, conjugation 
$\pi_f\mapsto {^\sigma\pi}_f$ defines a map
$$
\Pi_{\bf k}(\omega,K)\rightarrow \Pi_{^\sigma{\bf k}}({^\sigma\omega},K).\leqno(2.1)
$$
For any $\sigma\in{\rm Aut}({\Bbb C}/{\Bbb Q})$ we define a conjugation map on Hecke algebras:  
$$
\begin{array}{cccc}
\sigma:&{\cal H}_v(\omega_v,K_v)&\rightarrow&{\cal H}_v({^\sigma}\omega_v,K_v)\\
&\varphi&\mapsto&{^\sigma}\varphi,\\
\end{array}
$$
where $(\sigma(\varphi))(x)=\sigma(\varphi(x))$. Quite analogous we define a map
$$
\sigma:\,{\cal H}_f(\omega_f,K)\rightarrow{\cal H}_f({^\sigma}\omega_f,K)
$$
by sending $\otimes_{v\not\in S_\infty} \varphi_v$ to $\otimes_{v\not\in S_\infty} \sigma(\varphi_v)$. We note an easy property. Let $v\in V$. We denote by
$$
{\rm Ann}_\pi(v)=\{\varphi\in{\cal H}_v(\omega_v,K_v)\,(\mbox{resp.}\,\varphi\in {\cal H}_f(\omega_f,K)):\;\pi(\varphi)(v)=0\}
$$
the annihilator of $v$ in ${\cal H}_v(\omega_v,K_v)$ (resp. in ${\cal H}_f(\omega_f,K))$. Let $\varphi$ be in the Hecke algebra of 
$G({\Bbb A}_f)$ (resp. of $G(F_v)$); since
$$
^\sigma\pi(\varphi)(v)=\pi(\sigma^{-1}(\varphi))(v)
$$
for all $v\in {^\sigma W}=W$, we obtain for all $v$ 
$$
{\rm Ann}_{{^\sigma \pi}}\,v=\sigma({\rm Ann}_\pi \,v).\leqno(2.2)
$$

\medskip

{\bf 2.2. The Trace identity.}  We fix a compact open subgroup $K=\prod_{v\not\in S_\infty} K_v\le \prod_{v\not\in S_\infty} G({\cal O}_v)$ and a
character $\omega:\,F^*\backslash {\Bbb A}^*\rightarrow {\Bbb C}^*$. For any finite place $v$ we fix a compact open subgroup $K_{v,{\rm ell}}\le K_v$ and an element 
${\bf b}_{v,{\rm ell}}\in{\cal H}_v(\omega_v,K_{v,{\rm ell}})$ as follows. We denote by $G(F_v)_{\rm ell}$ the set of $F_v$-elliptic elements in 
$G(F_v)$. Lemma 7.2 (iii) in [J-L] implies that $G(F_v)_{\rm ell}\le G(F_v)$ and, hence, $K_v\cap G(F_v)_{\rm ell}\le K_v$ is an open subgroup. 
We select an (elliptic) element $\beta_v\in K_v\cap G(F_v)_{\rm ell}$; $K_v\cap G(F_v)_{\rm ell}$ then contains a neighbourhood 
$\beta_v K_{v,{\rm ell}}$ of $\beta_v$, where $K_{v,{\rm ell}}\le K_v\cap G(F_v)_{\rm ell}$ is an open group. Shrinking 
$K_{v,{\rm ell}}$ (intersect with principal congruence subgroups) we may assume that $zk=k'$ where $z\in Z(F_v)$ and 
$k,k'\in K_{v,{\rm ell}}\beta_v K_{v,{\rm ell}}$ implies that $z\in {\rm ker}\,\omega_v$ (to see this write 
$K_{v,{\rm ell}}\beta_v K_{v,{\rm ell}}=K'_v\beta_v$, where $K'_v=\beta_v K_{v,{\rm ell}}\beta_v^{-1}$ then is a subgroup of the
principal congruence subgroup of level $p^r$).
We define the function ${\bf b}_{v,{\rm ell}}$ on $G(F_v)$ by
$$
{\bf b}_{v,{\rm ell}}(x)=\left\{\begin{array}{cl}
0&x\not\in Z(F_v) K_{v,{\rm ell}}\beta_v K_{v,{\rm ell}} \\
\omega_v(z)^{-1}&x=zk\in Z(F_v)\cdot K_{v,{\rm ell}}\beta_v K_{v,{\rm ell}}.\\
\end{array}
\right.
$$
It is easy to see that ${\bf b}_{v,{\rm ell}}$ is well defined and  that
$$
{\bf b}_{v,{\rm ell}}\in{\cal H}_v(\omega_v,K_{v,{\rm ell}}).
$$
For any finite place $u$ we set 
$K_{\{u\}}=\prod_{v\not\in S_\infty,\,v\not=u} K_v\times K_{u,{\rm ell}}\,\le K$. Moreover, we set 
${\cal H}_{\{u\}}(\omega_f,K)=\otimes_{v\not|\infty,\,v\not=u} {\cal H}_v(\omega_v,K_v)\otimes{\bf b}_{u,{\rm ell}}$ and 
$$
{\cal H}_{\rm ell}(\omega_f,K)=\bigcup_{u\not\in S_\infty} {\cal H}_{\{u\}}(\omega_f,K)\;\subseteq {\cal H}_f(\omega_f,K).
$$ 
We define the conjugate of the weight ${\bf k}=(k_v)_{v\in S_\infty}$ by
$$
{^\sigma}{\bf k}=(k_{\sigma^{-1}(v)})_{v\in S_\infty}.
$$

\medskip

{\bf Proposition 2.1. }{\it Let $\sigma\in{\rm Aut}({\Bbb C}/{\Bbb Q})$. Assume that for all 
$\varphi\in{\cal H}_{\rm ell}(\omega_f,K)$ 
$$
\sigma\,{\rm tr}\,\varphi|_{L_0^2(\omega,{\bf k})}={\rm tr}\,{^\sigma}\varphi|_{L_0^2({^\sigma\omega},{^\sigma}{\bf k})}.\leqno(2.3)
$$
%
%
Then, conjugation by $\sigma\in{\rm Aut}({\Bbb C}/{\Bbb Q})$ defines a map
$$
\begin{array}{ccc}
\Pi_{\bf k}(\omega,K)&\rightarrow&\Pi_{^\sigma{\bf k}}({^\sigma\omega},K)\\
\pi&\mapsto&{^\sigma}\pi.\\
\end{array}
$$

}

\medskip

Thus, if $\pi$ is the finite part of a cuspidal automorphic representation of $G({\Bbb A})$, which has ${\rm SO}_2({\Bbb R})$-type given by
${\bf k}$ then ${^\sigma}\pi$ is the finite part of a cuspidal automorphic representation of $G({\Bbb A})$, which has 
${\rm SO}_2({\Bbb R})$-type given by ${\bf k}$. We deduce the Proposition from the following

\medskip

{\bf Lemma 2.2. } {\it Let $\sigma\in{\rm Aut}({\Bbb C}/{\Bbb Q})$ and fix an arbitrary finite place $u$. Assume that for all
$\varphi\in{\cal H}_{\{u\}}(\omega_f,K)$ we have
$$
\sigma\,{\rm tr}\,\varphi|_{L_0^2(\omega,{\bf k})}={\rm tr}\,{^\sigma}\varphi|_{L_0^2({^\sigma\omega},{^\sigma}{\bf k})}.
$$
Then, for any $\pi\in \Pi_{\bf k}(\omega,K)$ there is a representation 
$\pi'=\pi'_{\{u\}}\in \Pi_{^\sigma{\bf k}}({^\sigma\omega},K_{\{u\}})$ such that
$
{^\sigma\pi}_v\cong \pi'_{v}
$
for all finite places $v\not=u$. 


}

\medskip

Note that $\pi'$ is allowed to have smaller level than $\pi$.

{\it Proof of Lemma 2.2. } We set ${\cal H}_v={\cal H}_v(\omega_v,K_v)$, $v\not=u$. We write
$$
\Pi_{\bf k}(\omega,K_{\{u\}})=\{(\pi_1,V_1)\oplus\ldots\oplus(\pi_n,V_n)\}
$$
and 
$$
\Pi_{^\sigma{\bf k}}({^\sigma\omega},K_{\{u\}})=\{(\pi'_1,V_1')\oplus\ldots\oplus(\pi'_m,V_m')\}.
$$
We denote by $v_{i,v}\in V_{i,v}$ resp. $v_{i,v}'\in V_{i,v}'$ the essential vector, 
hence, $V_{i,v}^{K_v}={\cal H}_v v_{i,v}$ and $(V_{i,v}')^{K_v}={\cal H}_v v_{i,v}'$ for all $v\not=u$, because $\pi_i$ and $\pi_i'$
are irreducible. 
We set
$$
{\mathfrak a}_{i,v}={\rm Ann}_{\pi_{i,v}}(v_{i,v})\quad\mbox{and}\quad {\mathfrak a}_{i,v}'={\rm Ann}_{\pi_{i,v}'}(v_{i,v}')
$$
and obtain for all $\not=u$
$$
V_{i,v}^{K_{v}}\cong {\cal H}_v/{\mathfrak a}_{i,v}\qquad\mbox{and}\qquad 
(V_{i,v}')^{K_{v}}\cong {\cal H}_v/{\mathfrak a}_{i,v}'.\leqno(2.4)
$$ 
as ${\cal H}_v$-modules. Let $(\pi,V)\in\Pi_{\bf k}(\omega,K)$ and assume there 
is no $\pi'=\pi'_{\{u\}}\in\Pi_{^\sigma{\bf k}}({^\sigma\omega,K})$ such that $\pi'_v\cong{^\sigma\pi}_v$ for all $v\not=u$. 
Without loss of generality we may assume that $(\pi,V)=(\pi_1,V_1)$ and we set ${\mathfrak a}_v={\mathfrak a}_{1,v}$. 
Since the 
representations $\pi_1,\ldots,\pi_n$ are pairwise non-isomorphic, there is for any $i\ge 2$ a place $v_i\not=u$ such that 
$\pi_{v_i}\not\cong\pi_{i,v_i}$. Since the local representations $\pi_{i,v_i}$ are irreducible and $V_{i,v_i}^{K_{v_i}}\not=0$ we 
know that $\pi_{v_i}\cong\pi_{i,v_i}$ precisely if $\pi_{v_i}^{K_{v_i}}\cong \pi_{i,v_i}^{K_{v_i}}$ as ${\cal H}_{v_i}$-module. In 
particular, using equation (2.4) we obtain for all $i\ge 2$ that
$$
{\mathfrak a}_{v_i}\not={\mathfrak a}_{i,v_i}.
$$  
%
%
Since ${\mathfrak a}_{v_i},{\mathfrak a}_{i,v_i}\le {\cal H}_{v_i}$ are maximal ideals, the Chinese Remainder Theorem implies that 
for all $i\ge 2$ there is $\varphi_{i,v_i}\in{\cal H}_{v_i}$ such that 
$$
\varphi_{i,v_i}\equiv 1\pmod{{\mathfrak a}_{v_i}}\quad\mbox{and}\quad \varphi_{i,v_i}\equiv 0\pmod{{\mathfrak a}_{i,v_i}}.
$$ 
Hence, for all $i=2,\ldots,n$ we know that
$$
\pi_{v_i}(\varphi_{i,v_i})={\rm id}_{V_{v_i}^{K_{v_i}}}\quad\mbox{and}\quad \pi_{i,v_i}(\varphi_{i,v_i})=0_{V_{i,v_i}^{K_{v_i}}}.\leqno(2.5)
$$
On the other hand, our assumption implies that $^\sigma\pi$ is not isomorphic to any representation $\pi'_i$, hence for any $i\ge 1$ there is 
a place $v_i'\not=u$ such that $\pi_{v_i'}\not\cong {^{\sigma^{-1}}\pi}'_{i,v_i'}$. As above we see that for all 
$i=1,\ldots,m$ there is $\varphi_{i,v_i'}\in{\cal H}_{v_i'}$ such that
$$
\pi_{v_i'}(\varphi_{i,v_i'})={\rm id}_{V_{v_i'}^{K_{v_i'}}}\quad\mbox{and}\quad \big({^{\sigma^{-1}}\pi'_{i,v_i'}\big)
(\varphi}_{i,v_i'})=0_{{V_{i,v_i'}'}^{K_{v_i'}}}.\leqno(2.6)
$$
We then define the element 
$$
\varphi=\bigotimes_{i=2}^n\varphi_{i,v_i}\,\otimes \bigotimes_{i=1}^m \varphi_{i,v_i'}\otimes {\bf b}_{u,{\rm ell}}\in{\cal
H}_{\{u\}}(\omega_f,K).
$$
The choice of the local components of $\varphi$ implies

\begin{itemize}

\item  Equation (2.5) implies that $\pi_i(\varphi)=0$ for all $i\ge 2$

\item  Equation (2.6) implies that $\varphi_{i,v_i'}\in{\rm Ann}_{{^{\sigma^{-1}}}\pi_{i,v_i'}'}\,v_{i,v_i'}'$, hence, equation (2.2) shows that
${^\sigma}\varphi_{i,v_i'}\in \sigma({\rm Ann}_{^{\sigma^{-1}} \pi_{i,v_i'}'} v_{i,v_i'}')={\rm Ann}_{\pi_{i,v_i'}'}v_{i,v_i'}'$. We thus obtain
$$
\pi_i'({^\sigma}\varphi)=0.
$$ 
for all $i\ge 1$. 

\item Using equations (2.5) and (2.6) we see that
$$
\pi(\varphi)={\rm id}_{V^{K_{\{u\}}}}.
$$

\end{itemize}

Altogether we have shown that
$$
{\rm tr}\,\varphi|_{L_0^2({\bf k},\omega)^K}={\rm dim}\,V^{K_{\{u\}}}\quad\mbox{and}\quad{\rm tr}\,{^\sigma}\varphi|_{L_0^2({^\sigma\omega},{^\sigma}{\bf k})^K}=0.
$$
This contradicts our assumption and the Lemma is proven.


\medskip

{\it Proof of Proposition 2.1. } Let $\pi\in\Pi_{\bf k}(\omega,K)$. We choose pairwise distinct places $u,w$. Lemma 2.2 implies that there are 
representations $\pi_1'=\pi_{\{u\}}'$ resp. $\pi_2'=\pi_{\{w\}}'$ in $\Pi_{^\sigma{\bf k}}(^\sigma\omega,K_{\{u\}})$ resp. 
$\Pi_{^\sigma{\bf k}}(^\sigma\omega,K_{\{w\}})$ such that 

1. $\pi'_{\{u\},v}\cong {^\sigma\pi}_v$ for all $v\not=u$ and $\pi'_{\{w\},v}\cong{^\sigma\pi}_v$ for all $v\not=w$. 

2. ${\pi'_{\{u\}}}^{K_{\{u\}}}\not=0$ and ${\pi_{\{w\}}'}^{K_{\{w\}}}\not=0$.

The multiplicity $1$ Theorem implies that $\pi_{\{u\}}'\cong \pi_{\{w\}}'$. Hence, $\pi_{\{u\}}\cong \pi$ and
$\pi_{\{u\}}^{K}\not=0$. Thus, $\pi'_{\{u\}}\in \Pi_{^\sigma{\bf k}}(^\sigma,\omega,K)$ is the looked for representation.

\medskip

\section{Algebraicity of Distributions. }

{\bf 3.1. The Comparison. } We compare simple trace formulas for $\varphi\in{\cal H}_f(\omega_f,K)$ and $^\sigma\varphi\in {\cal H}_f({^\sigma\omega}_f,K)$ to show that equation (2.3) holds. We note that this will imply our main result on conjugation of Hilbert modular forms (cf. Theorem 3.3 and Corollary 3.4 below). 
As before, we fix a level $K=\prod_{v\not\in S_\infty} K_v\le \prod_{v\not\in S_\infty}G({\cal O}_v)$ 
and a character $\omega:\,F^*\backslash {\Bbb A}^*\rightarrow{\Bbb C}^*$ and we assume that condition (Alg) holds. 
Moreover, we set $\tilde{G}={\rm GL}_2/Z$ and we denote by 
$\tilde{G}_{F,{\rm ell}}$ the set of $F$-elliptic elements in $\tilde{G}(F)$. 

\medskip

{\bf Proposition 3.1. }{\it For any $\varphi\in{\cal H}_{\rm ell}(\omega_f,K)$ and any $\sigma\in{\rm Aut}({\Bbb C}/{\Bbb Q})$ we have
$$
\sigma\,{\rm tr}\,\varphi|_{L_0^2(\omega,{\bf k})}={\rm tr}\,{^\sigma}\varphi|_{L_0^2({^\sigma}\omega,{^\sigma}{\bf k})}.
$$

}

{\it Proof. } We define the following distributions on ${\cal H}_f(\omega_f,K)$: 
for any $\varphi\in{\cal H}_f(\omega_f,K)$ we set
$$
{\bf J}_{e,{\bf k},w}(\varphi)={\rm meas}\,\tilde{G}(F)\backslash \tilde{G}({\Bbb A})\;\varphi\otimes\hat{\varphi}_{{\bf k},w}(e)
$$
and
$$
{\bf J}_{{\rm ell},{\bf k},w}(\varphi)=\int_{\tilde{G}(F)\backslash \tilde{G}({\Bbb A})}\,
\sum_{\gamma\in \tilde{G}_{F,{\rm ell}}}\,\varphi\otimes\hat{\varphi}_{{\bf k},w}(x^{-1} \gamma x)\,dx,
$$ 
where $\hat{\varphi}_{{\bf k},w}$ is defined in section 1.4. Let $\varphi\in{\cal H}_{\rm ell}(\omega_f,K)$. 
For at least one finite place $u$ we know that $\varphi_u={\bf b}_{u,{\rm ell}}$. Since the local hyperbolic orbital integrals vanish for ${\bf b}_{u,{\rm ell}}$ as well as for $\varphi_{k_v}$, $v\in S_\infty$, the simple trace formula (cf. [G-J], Theorem 7.21, p. 245) yields
$$
{\rm tr}\,\tilde{\varphi}\otimes\varphi_{{\bf k},w}|_{L_d^2}(\omega)
={\bf J}_{e,{\bf k},w}(\varphi)+{\bf J}_{{\rm ell},{\bf k},w}(\varphi).
$$
The discrete spectrum decomposes $L_d^2(\omega,{\bf k})=L_d^2(\omega,{\bf k})\oplus L_{\rm res}^2(\omega,{\bf k})$, where $L_{\rm res}^2(\omega,{\bf k})$ vanishes unless ${\bf k}={\bf 2}=(2,\ldots,2)$ and $L_{\rm res}^2(\omega,{\bf 2})=\bigoplus_{\chi^2=\omega} \chi\circ{\rm det}$. 
Since $\varphi_{{\bf 2},w}$ has non-trivial ${\rm SO}_2({\Bbb R})$-type $\delta_{\bf 2}$ we know that the operator $\chi\circ{\rm det}(\varphi_{{\bf 2},w})$ vanishes, hence, ${\rm tr}\,\varphi\otimes\varphi_{{\bf 2},w}|_{L_{\rm res}^2(\omega,{\bf 2})}=0$. Together with Corollary 1.3 we obtain
$$
{\rm tr}\,{\varphi}|_{L_0^2(\omega,{\bf k})}
={\rm tr}\,{\varphi}\otimes\hat{\varphi}_{{\bf k},w}|_{L_0^2(\omega)}
={\rm tr}\,\varphi\otimes\varphi_{{\bf k},w}|_{L_d^2}(\omega)
={\bf J}_{e,{\bf k},w}(\varphi)+{\bf J}_{{\rm ell},{\bf k},w}(\varphi).
$$
Thus, the Propsition follows from the following Lemma, which completes the proof.

\medskip

{\bf Lemma 3.2. } {\it The distributions ${\bf J}_{e,{\bf k},w}$ and ${\bf J}_{{\rm ell},{\bf k},w}$ are algebraic, i.e.
$$
\sigma\,{\bf J}_{e,{\bf k},w}(\varphi)={\bf J}_{e,{\bf k},w}({{^\sigma}\varphi})\quad\mbox{and}\quad \sigma\,{\bf J}_{{\rm ell},{\bf k},w}(\varphi)={\bf J}_{{\rm ell},{\bf k},w}({{^\sigma}\varphi})
$$
for all $\varphi\in{\cal H}_f(\omega_f,K)$ and all $\sigma\in{\rm Aut}({\Bbb C}/{\Bbb Q})$.
}

\medskip


We will  give the proof of Lemma 3.2 in section 3.2. Proposition 2.1 and Proposition 3.1 imply our final result.

\medskip

{\bf Theorem 3.3. }{\it Let $\omega:\,F^*\backslash {\Bbb A}^*\rightarrow {\Bbb C}^*$  be an idele class character and let ${\bf k}=(k_v)_{v\in S_\infty}$ be a
weight such that (Alg) holds. Then conjugation by $\sigma\in{\rm Aut}({\Bbb C}/{\Bbb Q})$ defines a map
$$
\sigma:\,\Pi_{\bf k}(\omega,K)\rightarrow \Pi_{^\sigma{\bf k}}(^\sigma\omega,K).
$$
In different words, if $\pi_f$ is the finite part of a cuspidal automorphic representation of ${\rm GL}_2({\Bbb A})$ with central character $\omega$ and of lowest ${\rm SO}_{2,\infty}$-type $\delta_{\bf k}$, then, for any 
$\sigma\in{\rm Aut}({\Bbb C}/{\Bbb Q})$, ${^\sigma}\pi_f$ is the finite part 
of a cuspidal representation with central character ${^\sigma}\omega$ and lowest ${\rm SO}_{2,\infty}$-type $\delta_{{^\sigma}{\bf k}}$.

}

\medskip

{\bf Corollary 3.4. } {\it Any cuspidal representation $\pi\in\Pi_{\bf k}(\omega,K)$ is defined over a finite extension $E/{\Bbb Q}$.
}

{\it Proof. } We denote by ${\cal G}_{{\bf k},\omega}$ the set of all $\sigma\in{\cal G}$ satisfying ${^\sigma}{\bf k}={\bf k}$ and ${^\sigma}\omega=\omega$. Since the orbits ${\cal G}{\bf k}$ and ${\cal G}\omega$ are obviously finite
we know that $[{\cal G}:{\cal G}_{{\bf k},\omega}]<\infty$. Let $\pi\in\Pi_{\bf k}(\omega,K)$. Theorem 3.3 implies that ${^\sigma}\pi\in\bigcup_{\tau\in{\cal G}/{\cal G}_{{\bf k},\omega}}\Pi_{{^\tau}{\bf k}}({^\tau}\omega,K)$ for all $\sigma\in{\rm Aut}({\Bbb C}/{\Bbb Q})$. The union on the right hand side is a finite set, hence, the orbit ${\cal G}\pi$ is finite and we deduce that the stabilizer ${\cal G}_\pi$ of $\pi$ has finite index in ${\cal G}$. Since $\pi$ is defined over the subfield $E$ of ${\Bbb C}$, which is invariant under  ${\cal G}_\pi$ (cf. [W] or [C], Proposition 3.1) this proves the Corollary.

\medskip

{\bf 3.2. Proof of Lemma 3.2: Algebraicity of ${\bf J}_{{\rm ell},{\bf k},w}$ and  ${\bf J}_{e,{\bf k},w}$. } In this section we prove Lemma 3.2. The essential step will be to show that the local {\it archimedean} orbital integrals attached 
to $\hat{\varphi}_{{\bf k},w}$ are algebraic (cf. Lemma 3.4 below). To prove this, we use a  transfer of orbital integrals from the group $G={\rm GL}_2$ to 
the compact form $G'={\rm SU}_2({\Bbb C})$. This is the method of [Ca] in the case ${\rm PSL}_2$ and weight $k=2$, which we extend to the simply connected case
and arbitrary weight $k\ge 2$. We start by defining the local orbital integrals. We denote by $dx$ a Haar measure on ${\rm GL}_2({\Bbb R})$ or on  ${\rm SL}_2({\Bbb R})$ such that ${\rm meas}\, {\rm SO}_2({\Bbb R})=1$.
%
%
We set 
$$
T'=\{\Mat{\alpha}{}{}{\alpha^{-1}},\;\alpha\in{\Bbb C},\,\alpha\bar{\alpha}=1\}\cong S^1.
$$
Thus, $T'$ is a maximal torus in $G'$ and we choose Haar measures $d'g$ resp. $d't$ on $G'$ resp. on $T'$ such that 
${\rm meas}\,G'={\rm meas}\,T'=1$. Let ${\cal G}$ denote one of the groups $G={\rm SL}_2({\Bbb R})$ or $G'$
Let $t\in {\cal G}$ be a regular element and denote by ${\cal G}(\gamma)\le {\cal G}$ the centralizer of $\gamma$ in ${\cal G}$; 
for any function $\varphi$ on ${\cal G}$ we define the orbital integral
$$
{\cal O}_{\cal G}(\varphi,t)=\int_{{\cal G}(\gamma)\backslash {\cal G}} \varphi(x^{-1}t x)\,dx.
$$
Moreover, for any $x\in {\cal G}$ we set
$$
D(x)=\frac{|\alpha-\beta|_\infty^2}{|\alpha\beta|_\infty^2},
$$
where $\alpha,\beta$ are the eigenvalues of $x$. Finally, we denote by
$$
{\rm ch}_{k,w}:\,{\rm GL}_2({\Bbb C})\rightarrow{\Bbb C}
$$
the character of the finite dimensional irreducible representation $(L_{k,w},V_{k,w})$ of ${\rm GL}_2({\Bbb C})$ (cf. section 1.2). We 
denote by ${\rm ch}_k$ the restriction of ${\rm ch}_{k,w}$ to ${\rm SL}_2({\Bbb C})$, hence, ${\rm ch}_k$ is the character of $L_k$,
which is the restriction of $L_{k,w}$ to ${\rm SL}_2({\Bbb C})$.

\medskip

{\bf Lemma 3.5. }{\it Let $\gamma\in{\rm GL}_2({\Bbb R})$. If $\gamma$ is ${\Bbb R}$-hyperbolic or $\det \gamma<0$ then
$$
\int_{{\rm GL}_2({\Bbb R})(\gamma)\backslash{\rm GL}_2({\Bbb R})} \varphi_{k,w}(x^{-1}\gamma x)\,dx=0.
$$
If $\gamma$ is ${\Bbb R}$-elliptic and $\det \gamma>0$ then 
$$
\frac{1}{2}\,\int_{{\rm GL}_2({\Bbb R})(\gamma)\backslash{\rm GL}_2({\Bbb R})} \varphi_{k,w}(x^{-1}\gamma x)\,dx
=-{\rm ch}_{{k-2},w}(\gamma).
$$

}


%
%
%

{\it Proof. } Since $\varphi_{k,w}$ vanishes on all elements with negative determinant (cf. the proof of Lemma 1.1) and since the trace of $\varphi_{k,w}$ vanishes on principal series 
representations (cf. Corollary 1.2) the first claim is immediate. We therefore assume that $\gamma$ is ${\Bbb R}$-elliptic and $\det \gamma>0$. We first prove the statement of the Lemma in the case ${\rm SL}_2$ (cf. equation (3.5) below). As above we set 
$G={\rm SL}_2({\Bbb R})$ and $G'={\rm SU}_2({\Bbb C})$. For any $\gamma'\in G'$ there is an ${\Bbb R}$-elliptic element $\gamma$ in $G$ 
such that $\gamma$ and $\gamma'$ share the same characteristic polynomial. The assignment $\{\gamma'\}\mapsto\{\gamma\}$ defines a
bijection between the set of conjugacy classes in $G'$ and the set of ${\Bbb R}$-elliptic conjugacy classes in $G$. 
%
%
Any maximal torus in $G'$ is conjugate to $T'$. Similar, any non-split torus in $G$ is conjugate to $T={\rm SO}_2({\Bbb R})$. We denote 
by $(L_k',V_k')$ the restriction of the irreducible representation $(L_k,V_k)$ of ${\rm SL}_2({\Bbb C})$ of highest weight $k\gamma$ to 
${\rm SU}_2({\Bbb C})$ ($\gamma$ the fundamental root) and we denote by ${\rm ch}_{k}'$ the character of $L_k'$, hence, ${\rm ch}_{k}'={\rm ch}_{k}|_{G'}$. By Weyl's 
unitarian trick any irreducible representation of ${\rm SU}_2({\Bbb C})$ is of the form $(L_k',V_k')$ for some $k$. We set
$$
\chi_k=-\overline{{\rm ch}_{k-2}'},\quad k\ge 2
$$
(here, "$\bar{\quad}$" denotes the complex conjugation). ${\rm SU}_2({\Bbb C})$ is a compact group and the orthogonality relations imply for all $m\ge 0$
$$
{\rm tr}\,L_m'(\chi_k)=\left\{\begin{array}{ccc}
-{\rm meas}\,{\rm SU}_2({\Bbb C})&m=k-2\\
0&{\rm else.}\\
\end{array}\right.
$$
Comparing with Lemma 1.1 we obtain
$$
{\rm meas}\,{\rm SU}_2({\Bbb C})\;{\rm tr} D_m(\varphi_k)=-{\rm tr}\,L_{m-2}'(\chi_k)\leqno(3.1)
$$
for all $k,m\ge 2$. 
Since $\varphi_k$ vanishes on hyperbolic elements and since there is only one conjugacy class of non-split tori in $G$ we obtain using 
Weyl's integration formula (cf. [Kn], Proposition 5.27, p. 141) 
$$
{\rm tr}\, D_m(\varphi_k)=1/2\,\int_{T}\,D(t)\,{\rm ch}_{D_m}(t)\, {\cal O}_{G}(\varphi_k,t)\,dt.
$$
Quite analogous, Weyl's integration formula in the compact case (cf. [Kn], Theorem 4.45, p. 104) yields 
$$
{\rm tr}\, L_m'(\chi_k)=1/2\,\int_{T'}\,D(t)\,{\rm ch}_{m}'(t)\, {\cal O}_{G'}(\chi_k,t)\,d't.
$$
Hence, using (3.1) we deduce that
$$
{\rm meas}\,{\rm SU}_2({\Bbb C})\,\int_{T}\,D(t)\,{\rm ch}_{D_m}(t) {\cal O}_{G}(\varphi_k,t)\,dt
=-\int_{T'}\,D(t)\,{\rm ch}_{{m-2}}'(t)\, {\cal O}_{G'}(\chi_k,t)\,d't\leqno(3.2)
$$
for all $k,m\ge 2$. Lemma 1.1 implies that ${\rm ch}_{D_m}(t)=-{\rm ch}_{{m-2}}(t)$ for all $t\in {\rm SL}_2({\Bbb R})$. Moreover, if 
the conjugacy classes of $t\in{\rm SL}_2({\Bbb R})$ and $t'\in{\rm SU}_2({\Bbb C})$ correspond to each other, i.e. $t$ and $t'$ share 
the same characteristic polynomial, we know that ${\rm ch}_{m}(t)={\rm ch}_{m}(t')={\rm ch}_{m}'(t')$; 
hence, we obtain
$$
{\rm ch}_{D_m}(t)=-{\rm ch}_{{m-2}}'(t').\leqno(3.3)
$$
Together with equation (3.2) we obtain for all $m\ge 2$
$$
{\rm meas}\,{\rm SU}_2({\Bbb C})\,\int_{T}\,D(t)\,{\rm ch}_{{m-2}}'(t') {\cal O}_{G}(\varphi_k,t)\,dt
=\int_{T'}\,D(t)\,{\rm ch}_{{m-2}}'(t)\, {\cal O}_{G'}(\chi_k,t)\,d't.
$$
For any $t\in T$ there is $t'\in T'$ such that $t$ and $t'$ share the same characteristic polynomial and the assignment $t\mapsto t'$ 
defines a bijection $T\leftrightarrow T'$. Hence, Lemma 5.2.1 in [Ca] applied to $G'$ (note that there is only one conjugacy class of 
tori in $G'$ represented by $T'$) then implies for all $t\in T$ that
$$
{\cal O}_{G'}(\chi_k,t')={\rm meas}\,{\rm SU}_2({\Bbb C})\,{\cal O}_{G}(\varphi_k,t).\leqno(3.4)
$$
Since characters are class functions we obtain for all $t'\in T'$
$$
{\cal O}_{G'}(\chi_k,t')=\int_{T'\backslash G'}-\overline{{\rm ch}_{k-2}'}(x^{-1} t' x)\,dx
=- \frac{{\rm meas}\,{\rm SU}_2({\Bbb C})}{{\rm meas}\,T'}\,\overline{{\rm ch}_{k-2}'}(t').
$$
Hence, taking into account the normalization of measures we obtain using equation (3.4)
$$
{\cal O}_G(\varphi_k,t)=-\overline{{\rm ch}_{k-2}'}(t')=-\overline{{\rm ch}_{k-2}}(t).\leqno(3.5)
$$

Let now $\gamma\in{\rm GL}_2({\Bbb R})$ be ${\Bbb R}$-elliptic with positive determinant. We write $\gamma={\rm diag}(z,z)\gamma_0$ 
with $z\in{\Bbb R}^*_{>0}$ and $\gamma_0\in{\rm SL}_2({\Bbb R})$. Since $\gamma_0$ is ${\Bbb R}$-elliptic it is conjugate to an element 
$t\in {\rm SO}_2({\Bbb R})$. Since ${\rm GL}_2({\Bbb R})(\gamma)={\rm SL}_2({\Bbb R})^\pm(\gamma) Z_2^0$ we obtain
$$
\int_{{\rm GL}_2({\Bbb R})(\gamma)\backslash{\rm GL}_2({\Bbb R})} \varphi_{k,w}(x^{-1}\gamma x)\,dx
=\int_{{\rm SL}_2({\Bbb R})^\pm(\gamma)\backslash{\rm SL}_2({\Bbb R})^\pm} \varphi_{k,w}(x^{-1}\gamma x)\,dx.
$$

Since ${\rm SL}_2({\Bbb R})^\pm(t)={\rm SO}_2({\Bbb R})$ for all regular $t\in {\rm SO}_2({\Bbb R})$ 
%
and since $\gamma_0$ is conjugate to an element $t\in {\rm SO}_2({\Bbb R})$ we deduce that
${\rm SL}_2({\Bbb R})^\pm(\gamma)$ is conjugate to ${\rm SL}_2({\Bbb R})^\pm(t)={\rm SL}_2({\Bbb R})(t)$ and we obtain
$$
{\rm SL}_2({\Bbb R})^\pm(\gamma)={\rm SL}_2({\Bbb R})(\gamma).
$$
We deduce that
$$
{\rm SL}_2({\Bbb R})^\pm(\gamma) \backslash {\rm SL}_2({\Bbb R})^\pm
={\rm SL}_2({\Bbb R})(\gamma) \backslash {\rm SL}_2({\Bbb R}) \cup j {\rm SL}_2({\Bbb R})(j^{-1}\gamma j) \backslash {\rm SL}_2({\Bbb R}).
$$
Using the this decomposition we obtain
\begin{eqnarray*}
&&\int_{{\rm SL}_2({\Bbb R})^\pm(\gamma)\backslash{\rm SL}_2({\Bbb R})^\pm} \varphi_{k,w}(x^{-1}\gamma x)\,dx\\
&&=\int_{{\rm SL}_2({\Bbb R})(\gamma)\backslash{\rm SL}_2({\Bbb R})} \varphi_{k,w}(x^{-1}\gamma x)\,dx
+\int_{{\rm SL}_2({\Bbb R})(j^{-1}\gamma j)\backslash{\rm SL}_2({\Bbb R})} \varphi_{k,w}(x^{-1}j^{-1}\gamma j x)\,dx\\
&&=z^w \, \int_{{\rm SL}_2({\Bbb R})(\gamma)\backslash{\rm SL}_2({\Bbb R})} \varphi_{k}(x^{-1}\gamma_0 x)\,dx
+z^w \, \int_{{\rm SL}_2({\Bbb R})(j^{-1}\gamma j)\backslash{\rm SL}_2({\Bbb R})} \varphi_{k}(x^{-1}j^{-1}\gamma_0 j x)\,dx\\
\end{eqnarray*}
(note that $\varphi_{k,w}$ and $\varphi_k$ coincide on ${\rm SL}_2({\Bbb R})$). Using equation (3.5) this yields
\begin{eqnarray*}
\int_{{\rm SL}_2({\Bbb R})^\pm(\gamma)\backslash{\rm SL}_2({\Bbb R})^\pm} \varphi_{k,w}(x^{-1}\gamma x)\,dx 
&=&-z^w \overline{{\rm ch}_{k-2}}(\gamma_0)-z^w \overline{{\rm ch}_{k-2}}(j^{-1} \gamma_0 j)\\
&=&-\overline{{\rm ch}_{k-2,w}}(\gamma)-\overline{{\rm ch}_{k-2,w}}(j^{-1} \gamma j).\\
&=&-2\overline{{\rm ch}_{k-2,w}}(\gamma).\\
\end{eqnarray*}
(note that ${\rm ch}_{k-2,w}$ is constant on conjugacy classes). Since $L_{k,w}$ is defined over ${\Bbb Q}$ and $\gamma\in{\rm GL}_2({\Bbb R})$ we know that $\overline{{\rm
ch}_{k-2,w}}(\gamma)={\rm ch}_{k-2,w}(\overline{\gamma})={\rm ch}_{k-2,w}(\gamma)$. This completes the proof of the Lemma.

\medskip

We set $(L_{\bf k},V_{\bf k})=\otimes_{v\in S_\infty} (L_{k_v},V_{k_v})$, hence, $L_{\bf k}$ is a representation of ${\rm GL}_2({\Bbb
R})^{|S_\infty|}$ and we denote by 
$$
{\rm ch}_{{\bf k}}=\prod_{v\in S_\infty} {\rm ch}_{{k_v}}
$$
its the character. Using the identifications from section 1.1, we say that $\gamma\in G(F)$ is totally elliptic if $i_v(\gamma)\in G(F_v)$ is ${\Bbb R}$-elliptic for all $v\in S_\infty$ 
and totally positive if $\det i_v(\gamma)>0$ for all $v\in S_\infty$
We embed 
$$
i_\infty:\,G(F)\hookrightarrow G_\infty:=\prod_{v\in S_\infty} G(F_v),
$$ 
by $\gamma\mapsto(i_v(\gamma))_{v\in S_\infty}$. An immediate consequence of  Lemma 3.5 is

\medskip

{\bf Corollary 3.6. }{\it 1.) Let $\gamma\in G(F)$. If $\gamma$ is totally elliptic and totally positive then
$$
I_{{\bf k},w}(\gamma):=\int_{G_\infty(\gamma)/G_\infty}\hat{\varphi}_{{\bf k},w}(x^{-1} \gamma x)\,dx=-2\,{\rm ch}_{L_{{\bf k}-2},w}(\gamma),
$$
where ${\bf k}-2=(k_v-2)_{v\in S_\infty}$. Otherwise, the above integral vanishes. 

\medskip

2.) In particular, $I_{{\bf k},w}(\gamma)$ is algebraic, i.e. for all $\sigma\in{\rm Aut}({\Bbb C}/{\Bbb Q})$ and all $\gamma\in G(F)$ we have
$$
\sigma I_{{\bf k},w}(\gamma)=I_{^\sigma{\bf k},w}(\gamma).
$$

}

{\it Proof. } 1.) We compute using equation (1.1) and the Remark in section 1.1 
\begin{eqnarray*}
\int_{G_\infty(i_\infty(\gamma))\backslash G_\infty}\hat{\varphi}_{{\bf k},w}(x^{-1}i_\infty(\gamma)x)\,dx
&=&\prod_{v\in S_\infty}\int_{G(F_v)(i_v(\gamma))\backslash G(F_v)}\hat{\varphi}_{{k_v},w,v}(x^{-1}i_v(\gamma) x)\,dx\\
&=&\prod_{\tau\in{\cal G}}\int_{G(F_{\tau v_0})(i_{\tau v_0}(\gamma))\backslash G(F_{\tau v_0})} \hat{\varphi}_{{k_{\tau v_0}},w,\tau v_0}(x^{-1}i_{\tau
v_0}(\gamma) x)\,dx\qquad (v=\tau v_0)\\
&=&\prod_{\tau\in {\cal G}}\int_{{\rm GL}_2({\Bbb R})(\tau(\gamma))\backslash {\rm GL}_2({\Bbb R})}{\varphi}_{k_{\tau v_0},w}(x^{-1}\tau(\gamma)
x)\,dx\qquad \mbox{(Remark in sec. 1.1)}\\
&=&\prod_{\tau\in {\cal G}} {\rm ch}_{{k_{\tau v_0}},w}(\tau(\gamma))\qquad \mbox{(Lemma 3.5)}\\
&=&\prod_{\tau\in {\cal G}} {\rm ch}_{{k_{\tau v_0}},w}(i_{\tau v_0}(\gamma))\qquad \mbox{(eq. (1.1))}\\
&=&{\rm ch}_{{\bf k},w}(\gamma).\\
\end{eqnarray*}

\medskip

2.) Since the representation $L_{k,w}$ is defined over ${\Bbb Q}$, i.e. $\sigma\,L_{k,w}(g)=L_{k,w}({^\sigma g})$ for $g\in{\rm GL}_2({\Bbb C})$ 
we obtain $\sigma\,{\rm tr}\,L_{{\bf k},w}(\gamma)={\rm tr}\, L_{^\sigma{\bf k},w}(\gamma)$
%
%
%
Together with part 1.) this implies the claim and the Corollary is proven.

\medskip

We give the {\it Proof of Lemma 3.2. } 

{\it Algebraicity of ${\bf J}_{{\rm ell},{\bf k},w}$. } We write $\tilde{G}_\infty(\gamma)$ resp. $\Gamma(\gamma)$ for the centralizer of $\gamma$ in 
$\tilde{G}_\infty=\prod_{v\in S_\infty} \tilde{G}(F_v)$ resp. in $\Gamma$. We also denote by $\{\tilde{G}_{F,{\rm ell}}\}_\Gamma$ the set of $\Gamma$-conjugacy classes in $\tilde{G}_{F,{\rm ell}}$ and by $\{\gamma\}_\Gamma$ the $\Gamma$-conjugacy class of $\gamma\in \tilde{G}_{F,{\rm ell}}$. Strong approximation for ${\rm GL}_2$ implies that there is a finite set 
${\cal V}\subset\tilde{G}({\Bbb A}_f)$ such that
$$
\tilde{G}({\Bbb A})=\dot{\bigcup}_{\xi\in {\cal V}} \tilde{G}(F)\xi \tilde{G}_\infty \tilde{K}.
$$
We denote by $\tilde{K}$ the image of $K$ in $\tilde{G}({\Bbb A}_f)$ and by 
$\Gamma_\xi=\{\gamma\in \tilde{G}(F):\,\xi^{-1}\gamma\xi\in \tilde{K}\}$, $\xi\in{\cal V}$, the arithmetic subgroup of $\tilde{G}_\infty$ corresponding to $\tilde{K}$. Since the assignment
$$
x\mapsto \sum_{\gamma\in\tilde{G}_{F,{\rm ell}}}\, \varphi\otimes\hat{\varphi}_{{\bf k},w}(x^{-1}\gamma x)
$$ 
defines a function on $\tilde{G}({\Bbb A})$, which is left $\tilde{G}(F)$ 
and right $\tilde{K}$-invariant 
we obtain



\begin{eqnarray*}
&&\int_{\tilde{G}(F)\backslash \tilde{G}({\Bbb A})}\,\sum_{\gamma\in\tilde{G}_{F,{\rm ell}}}\,\varphi\otimes\hat{\varphi}_{{\bf k},w}(x^{-1} \gamma x)\,dx \\
&&\qquad={\rm vol}(K)\;\sum_\xi\,\int_{\Gamma_\xi\backslash \tilde{G}_\infty} \sum_{\gamma\in \tilde{G}_{F,{\rm ell}}} \varphi(\xi^{-1}\gamma\xi)
\,\hat{\varphi}_{{\bf k},w}(x^{-1} \gamma x)\,dx\\
&&\qquad={\rm vol}(K)\,\sum_\xi\,\sum_{\{\gamma\}_{\Gamma_\xi}\in \{\tilde{G}_{F,{\rm ell}}\}_{\Gamma_\xi}}\;\varphi(\xi^{-1}\gamma\xi)\;\int_{\Gamma_\xi\backslash
\tilde{G}_\infty} \sum_{\tau\in\Gamma_\xi(\gamma)\backslash \Gamma_\xi} \hat{\varphi}_{{\bf k},w}(x^{-1} \tau^{-1}\gamma\tau x)\,dx\\
&&\qquad={\rm vol}(K)\,\sum_\xi\,\sum_{\{\gamma\}_{\Gamma_\xi}\in \{\tilde{G}_{F,{\rm
ell}}\}_{\Gamma_\xi}}\;\varphi(\xi^{-1}\gamma\xi)\;\int_{\Gamma_\xi(\gamma)\backslash \tilde{G}_\infty} \hat{\varphi}_{{\bf k},w}(x^{-1} \gamma x)\,dx\\
&&\qquad={\rm vol}(K)\,\sum_\xi\,\sum_{\{\gamma\}_{\Gamma_\xi}\in \{\tilde{G}_{F,{\rm ell}}\}_{\Gamma_\xi}} {\rm
meas}\,\tilde{G}_\infty(\gamma)/\Gamma_\xi(\gamma)\,\;\varphi(\xi^{-1}\gamma\xi)\;\int_{\tilde{G}_\infty(\gamma)\backslash \tilde{G}_\infty} \hat{\varphi}_{{\bf
k},w}(x^{-1} \gamma x)\,dx.\\
\end{eqnarray*}
Since $\tilde{G}_\infty(\gamma)=\prod_{v\in S_\infty} \tilde{G}(F_v)(\gamma)$ is conjugate to ${\rm PSO}_{2,\infty}=\prod_{v\in S_\infty} {\rm PSO}_2({\Bbb R})$ 
and taking into account the natural bijection $\tilde{G}_\infty(\gamma)\backslash \tilde{G}_\infty\cong G_\infty(\gamma)\backslash G_\infty$ we obtain
\begin{eqnarray*}
{\bf J}_{{\rm ell},{\bf k},w}(\varphi)=
{\rm vol}(K\times {\rm PSO}_{2,\infty})\,\sum_\xi\,\sum_{\{\gamma\}_{\Gamma_\xi}\in \{\tilde{G}_{F,{\rm ell}}\}_{\Gamma_\xi}} |\Gamma_\xi(\gamma)|^{-1}
\;\varphi(\xi^{-1}\gamma\xi)\;I_{{\bf k},w}(\gamma).\\
\end{eqnarray*}
Using Corollary 3.6 the above equation immediately implies that ${\bf J}_{{\rm ell},{\bf k},w}$ is algebraic.


\medskip

{\it Algebraicity of ${\bf J}_{e,{\bf k},w}$. } The Plancherel Theorem for ${\rm SL}_2({\Bbb R})$ (cf. [Kn], Theorem 11.6, p. 401) together
with Lemma 1.1 yields
$$
\varphi_{k_v}(e)=\frac{1}{2\pi}\,(k_v-1)\,{\rm trace}\,D_{k_v}(\varphi_{k_v})=\frac{k_v-1}{2\pi},
$$
hence,
$$
\hat{\varphi}_{{\bf k},w}(e)=\left(\frac{1}{2\pi}\right)^{[F:{\Bbb Q}]}\,\prod_{v\in S_\infty}(k_v-1).
$$
On the other hand, we set $\tilde{K}_{\rm max}=\prod_{v\not\in S_\infty} \tilde{G}({\cal O}_v)$ and we normalize the Haar measure $dg_f$ on $\tilde{G}({\Bbb A}_f)$ 
such that ${\rm meas}\,\tilde{K}_{\rm max}\times{\rm PSO}_{2,\infty}=1$. Since
${\rm meas}\,(G(F)\backslash G({\Bbb A}))={\rm meas}\,(G(F)\backslash G({\Bbb A})/{\rm PSO}_{2,\infty}K_{\rm max})$
and since $G(F)\backslash G({\Bbb A})/{\rm PSO}_{2,\infty}K_{\rm max}$ is a finite union of spaces of the form $\Gamma\backslash{\cal H}^d$, which have volume $c (2\pi)^{[F:{\Bbb Q}]}$, where $c\in{\Bbb Q}$ (cf. [Fr], section 2.5, Proposition 5.1 and [Hi], section 2.7, Corollary 1, p. 71), we obtain
$$
\varphi\otimes\hat{\varphi}_{{\bf k},w}(e)=\varphi(e)\,c\, \prod_{v\in S_\infty}(k_v-1).
$$
Since $c\in{\Bbb Q}$ this clearly implies that ${\bf J}_{e,{\bf k},w}$ is algebraic and completes the proof of Lemma 3.2.

\medskip

\medskip

{\bf 3.3. Variants of the proof. }

$\bullet$ {\it The Kazhdan-Flicker trace formula. } We denote by $\Pi_{\bf k}^{\rm cusp}(\omega,K)$ the set of all representations 
occuring in $L_0^2(\omega,{\bf k})^K$, which are cuspidal at at least one finite place. We fix a finite place $u$ and a cuspidal ireducible representation 
$\rho$ of $G_u$ and we denote by ${\cal H}_f(\rho,\omega_f,K)$ the set of all $\varphi=\otimes_{v\not\in S_\infty} \varphi_v\in{\cal H}_f(\omega_f,K)$ such that 
$\varphi_u$ equals a matrix coefficient of $\rho$. Hence, $\varphi$ is discrete and cuspidal in the sense of [F-K], p. 191. The assignment 
$\varphi\mapsto^\sigma\varphi$ defines a map
$$
{\cal H}_f(\rho,\omega,K)\rightarrow {\cal H}_f({^\sigma\rho},^\sigma\omega,K).
$$

\medskip

{\bf Proposition 3.7. }{\it Assume that for all finite places $u$, all cuspidal representations $\rho$ of $G_u$, all $\varphi\in {\cal H}_f(\rho,\omega,K)$ and all $\sigma\in{\rm Aut}({\Bbb C}/{\Bbb Q})$
$$
\sigma\,{\rm tr}\,\varphi|_{L_0^2(\omega,{\bf k})^K}={\rm tr}\,{^\sigma}\varphi|_{L_0^2(^\sigma\omega,^\sigma{\bf k})^K}.\leqno(3.7)
$$
Then conjugation defines a map
$$
\Pi_{\bf k}^{\rm cusp}(\omega,K)\rightarrow \Pi_{^\sigma{\bf k}}^{\rm cusp}(^\sigma\omega,K).
$$
More precisely, denote by $\Pi_{\bf k}^\rho(\omega,K)$ the set of all representations $\pi_f$ occuring in $L_0^2(\omega,{\bf k})^K$ and satisfying 
$\pi_u\cong\rho$; then conjugation defines a map
$$
\Pi_{\bf k}^\rho(\omega,K)\rightarrow \Pi_{^\sigma{\bf k}}^{{^\sigma}\rho}(^\sigma\omega,K).
$$
}

We sketch the {\it proof} of the second claim following the argument in the proof of Lemma 2.2. Let $\pi_f\in\Pi_{\bf k}^\rho(\omega,K)$ and assume that there is no $\pi'\in\Pi_{^\sigma{\bf k}}^{{^\sigma}\rho}(^\sigma\omega,K)$ such that ${^\sigma\pi}\cong\pi'$. Hence, for any representation 
$\pi'\in \Pi_{^\sigma{\bf k}}^{^\sigma\rho}(^\sigma\omega,K)$ there is a finite place $v\not=u$ such that ${^{\sigma^{-1}}}\pi_u'\not\cong \pi_u$. Moreover, for any 
representation $\xi\in \Pi_{\bf k}^\rho(\omega,K)$, $\xi\not\cong \pi$, there is a finite place $v\not=u$ such that $\pi_u\not\cong \xi_u$. As in the proof of Lemma 2.2 
we then construct a function $\varphi=\otimes_{v\in S_\infty}\varphi_v$ such that 
- $\varphi_u$ equals a matrix coefficient of $\rho$
- ${\rm tr}\,\xi(\varphi)=0$ for all $\xi\in \Pi_{\bf k}^{\rho}(\omega,K)$, which are not isomorphic to $\pi$ 
- ${\rm tr}\,\pi'(\varphi)=0$ for all $\pi'\in\Pi_{^\sigma{\bf k}}^{^\sigma\rho}(^\sigma\omega,K)$
- ${\rm tr}\,\pi(\varphi)\not=0$. 
This then contradicts our assumption (3.7) and proves the Proposition.

\medskip

Since $\varphi\in {\cal H}_f(\rho,\omega_f,K)$ is discrete and cuspidal in the sense of [F-K], p.191, we can apply the Kazhdan-Flicker simple trace formula, 
which yields
$$
{\rm tr}\,\varphi\otimes\hat{\varphi}_{{\bf k},w}|_{L_0^2(\omega,{\bf k})}={\bf J}_{{\rm ell},{\bf k},w}(\varphi).
$$
Thus, the distribution ${\bf J}_{e,{\bf k},w}$ does not appear and we do not have to use the Plancherel Theorem. Moreover, the proof of the 
Kazhdan-Flicker simple trace formula is elementary compared to the proof of the Selberg trace formula.


$\bullet$ {\it The simple trace formula. } If $F/{\Bbb Q}$ is a totally real field of degree $\ge 2$ then we can apply the simple 
trace formula to any element $\varphi\otimes\hat{\varphi}_{{\bf k},w}$ with $\varphi\in{\cal H}_f(\omega_f,K)$ arbitrary, because the hyperbolic orbital integrals 
vanish for $\varphi_{k_v,w}$, $v\in S_\infty$ and $|S_\infty|\ge 2$, i.e. we do not have to assume that $\varphi$ is elliptic at one place. In particular, unlike in the proof of Proposition 2.1 (and Lemma 2.2) we do not need the Multiplicity 1 Theorem (we need not fix an auxiliary place $u$).

\medskip

$\bullet$ {\it The Selberg trace formula. } We can also try to establish (0.2) by using the (full) Selberg trace formula. In this case we do not
have to use Multiplicity 1 even in the case $F={\Bbb Q}$. On the other hand the geometric side is more complicated: we have to prove
algebraicity for distributions attached to unipotent and hyperbolic conjugacy classes and also for distributions attached to Eisenstein series.

\bigskip

{\Large\bf References}

[Ca] Casselman, W., The Hasse-Weil $\zeta$-function of some moduli varieties of dimension greater than one. Proceedings of Symposia in Pure Mathematics, Vol. {\bf 33} (1979), part 2, pp. 141 - 163 

[C] Clozel, L., Motifs et Formes Automorphes. Application du Principe de Fonctorialite. In: Automorphic Forms, Shimura Varieties and $L$-functions, Perspectives in Math {\bf 10} Academic Press 1990.


[F-K] Flicker, Y., Kazhdan, D., A simple trace formula. Journal d'Analyse Mathematique {\bf 50} (1998), pp. 189  200

[Fr] Freitag, E., Hilbert modular forms, Springer 1990

[G] Garrett, P., Holomorphic modular forms, Wadsworth and Brooks/Cole Advanced Books, 1990
 
[G-J] Gelbart, S., Jacquet, H., Forms of ${\rm GL}_2$ from the analytic point of view, Proc. Symp. in Pure Math., {\bf 33} (1979) 213- 251
 
[J-L] Jacquet, H., Langlands, R., Automorphic forms on ${\rm GL}_2$, LNM {\bf 114}, Springer 1970
 
[H] Hida, H., Elementary Theory of $L$-functions and Eisenstein series, Cambridge University press, 1993
 
[Kn] Knapp, A., Representation theory of semi simple Lie groups. Princeton Landmarks in Mathematics, Princeton University Press, 1986

[Sh], Shimura, G., On some arithmetic properties of modular forms of one and several variables, Ann of. Math (2) {\bf 102} (1975), no. 3, 491 - 515

[W] Waldspurger, J.-L., Quelques propri{\'e}t{\'e}s arithm{\'e}tiques  de certaines formes automorphes sur ${\rm GL}(2)$, Comp. Math. {\bf 54} (1985), 121 - 171

\newpage

\end{document}